\documentclass[11pt]{article}
\usepackage{amsmath}
\usepackage{amsthm}
\usepackage{amscd}
\usepackage{amssymb}

\title{A strong Schottky Lemma for nonpositively curved singular
spaces\\
{\em \small To John Stallings on his sixty-fifth birthday}} 
\author{
Roger C. Alperin, 
Benson Farb\thanks{Supported
in part by the NSF and by the Alfred P.\ Sloan Foundation.}, 
and 
Guennadi A. Noskov\thanks{Supported by GIF-grant G-454-213.06/95}}

\newtheorem{theorem}{Theorem}[section]
\newtheorem{proposition}[theorem]{Proposition}
\newtheorem{lemma}[theorem]{Lemma}
\newtheorem{corollary}[theorem]{Corollary}

\newtheorem{exmple}[theorem]{Example}

\def\proof{{\bf {\medskip}{\noindent}Proof. }}

\def\endproof{$\diamond$ \bigskip}

\def\title{\em}

\newcommand\Z{\mbox{\bf Z}}
\newcommand\Q{\mbox{\bf Q}}

\newcommand\GL{{\rm GL}}

\newcommand\s{\sigma}
\newcommand\val{\nu}%the valuation
\newcommand\Isom{\mbox{Isom}}
\newcommand\CAT{\mbox{CAT}}

\DeclareMathOperator{\Lk}{Lk}

\begin{document}

\maketitle

\section{Introduction}

The classical Schottky Lemma (due to Poincare, Klein, Schottky) gives a 
criterion for a pair of isometries $g,h$ of hyperbolic space to have 
powers $g^m,h^n$ which generate a free group.  This criterion was
generalized by Tits to pairs of elements in linear groups in his proof 
of the Tits alternative.

In this paper we give a criterion (Theorem \ref{theorem:main}) for pairs
of isometries of a nonpositively curved metric space (in the sense of
Alexandrov) to generate a free group {\em without having to take
powers}.  This criterion holds only in singular spaces, for example in
Euclidean buildings; in fact our criterion takes a particularly simple
form in that case (Corollary \ref{theorem:buildings}).

The original motivation for our criterion was to prove that the four
dimensional Burau representation is faithful. While linearity of braid
groups is now known, this question is still open; it is well-known to be
related to detecting the unknot with the Jones polynomial .  It was
shown in \cite{Bir} that the faithfulness question is equivalent to
proving that a specific pair of elements in $\GL_3(\Z[x,x^-1])$ generate
a free group.  In \S\ref{section:examples} we show that these elements
are the image of the free group $F_2$ under a representation
$\rho_{0,0}$ lying in a $2$-parameter family of representations
$$\rho_{\alpha,\beta}:F_2\rightarrow
\GL_3(\Q(x))$$
\noindent
for $\alpha,\beta\in\Q$.  Our criterion applies to give
faithfullness of $\rho_{\alpha,\beta}$ for all $\alpha,\beta$ except, 
unfortunately, for $\alpha,\beta \in\{0,1\}$.  

\bigskip
\noindent
{\bf The criterion. }
For definitions of the terms we use we refer the reader to
\S\ref{section:cat0}.  
A complete, nonpositively curved (i.e.\ $\CAT(0)$) metric space $X$ is
said to have {\em no fake angles} if there is no pair of geodesics
which issue from the same point, have zero angle at that point, yet are
disjoint except at that point.  This is a very weak condition, and is
satisfied for example by any piecewise Euclidean $\CAT(0)$ simplicial
complex with finitely many isometry types of cells (Proposition
\ref{propsec}).

\begin{theorem}[Strong Schottky]
\label{theorem:main}
Let $X$ be a complete CAT(0)-space with no fake zero angles.  
Let $g_1,g_2$ be the axial isometries of $X$ with axes 
$A_1,A_2$ respectively.  
Assume that either one of the following holds:
\begin{enumerate}
\item $S=A_1\cap A_2$ is a bounded segment. 
 Let $s_-,s_+$ be its  endpoints and let   $A_1^\pm, A_2^\pm$ 
be the  infinite rays of $A_1,A_2$ with the origins $s_-,s_+$.
 Assume that the angles between  $A_1^+,A_2^+$
and between $A_1^-,A_2^-$ are equal to $\pi$, and that 
the translation lengths of $g_1,g_2$ are strictly
greater than the length of  $S.$

\item $A_1,A_2$ are disjoint, and  there is a geodesic
$B$  connecting $a_1\in A_1$ to $a_2\in A_2$ with all 
four angles between $B$ and  $A_1,A_2$ at $a_1,a_2$ equal to  $\pi$.
\end{enumerate}
Then the group generated by $g_1,g_2$ is free.
\end{theorem}

\bigskip
\noindent
{\bf The case of buildings. }Let $X$ be a thick affine building. For
example, let $X=\cal T_\val$ be the building associated to $\GL_n(K), 
n\ge 2$ where $K$ is a discretely-valued field with valuation
$\val$.  The link $link(x)$ of a point $x\in X$ is a spherical building.
Two chambers in the link are called {\em opposite} if their distance is
maximal in the link.  Recall that in any apartment $A$ of $X$, any
chamber in $A$ has a unique opposite chamber in $A$.
Recall that for any apartment $A$, any chamber of $A\cap link(x)$
has a unique opposite chamber in
$A\cap link(x)$.

An isometry $f\in \Isom(X)$ of hyperbolic type is said to be {\em
generic} if none of its (parallel) axes are contained in any wall of any
apartment of $X$.  Note that $f$ is generic if and only if it has a
unique invariant apartment $P_f$. The axis of $f$ is contained in the union of two sectors 
$P_f^{+}, P_f^{-}$ of $P_f$ which are invariant by translation by $f$, respectively $f^{-1}$. A generic
isometry $f$ determines, for any fixed choice of basepoint $x\in P_f$, a
pair of chambers in $link(x)$.  We say that generic $f,g\in
\Isom(X)$ are {\em opposite} if $P_f\cap P_g=x\in X$ and if each of the
chambers determined by $f$ is opposite in $link(x)$ to each of the chambers
determined by $g$.

As a straightforward consequence of Theorem \ref{theorem:main} we have 
the following.

\begin{corollary}[Strong Schottky for buildings]
\label{theorem:buildings}
Let $X$ be a (thick) Euclidean building, and let $f,g\in \Isom(X)$.  If 
$f,g$ are opposite then they generate a free subgroup of $\Isom(X)$.
\end{corollary}

In \S\ref{section:examples} we will apply 
Corollary \ref{theorem:buildings} to a $2$-parameter family of pairs of 
elements in the affine building associated to $\GL_3(\Q(x))$.  

\section{$\CAT(0)$ preliminaries}
\label{section:cat0}

\subsection{Definitions}

Let $X$ be a geodesic metric space.  The {\em comparison triangle} for a
geodesic triangle $\Delta$ in $X$ is the Euclidean triangle $\Delta'$
with the same side lengths as $\Delta$.  We say that $X$ is $CAT(0)$, 
or {\em nonpositively curved}, if for any geodesic triangle $\Delta$ in
$X$ and any two points $x,y$ on $\Delta,$ the distance between $x$ and
$y$ in $X$ is less than or equal to the Euclidean distance between the
corresponding points $x',y'$ on the comparison triangle $\Delta'$ ( 
\cite{BH}, II.1.1).  The $\CAT(0)$ condition implies the 
uniqueness of geodesics and the geodesicity of local geodesics.  In the
following we assume that $X$ is a complete $CAT(0)$ space.

Let $\varepsilon>0$ and let
$\sigma_1,\sigma_2:[0,\varepsilon]\rightarrow X$ be two unit speed
geodesics with $\sigma_1(0)=\sigma_2(0)=:x$.  For
$s,t\in(0,\varepsilon)$ let $\gamma(s,t)$ be the angle at $x'$ of the
comparison triangle $x'\sigma_1(s)'\sigma_1(t)'.$ Then $\gamma(s,t)$ is
monotonically decreasing as $s,t$ decrease and hence $
\angle(\sigma_1,\sigma_2):=\lim_{s,t\rightarrow 0}\gamma(s,t)
$ exists and is called the {\em angle} subtended by $\sigma_1$ and
$\sigma_2$ (\cite{Ba}, I.(3.8)).  A metric expression for the angle
is given by the ``cosine theorem''(\cite{Ba},
I.(3.11) $$
\cos(\angle(\sigma_1,\sigma_2))=
\lim_{s,t\rightarrow 0}
\frac
{s^2+t^2-d^2(\sigma_1(s),\sigma_2(t))} {2st} .$$ 

  It follows from this formula that if the angle is strictly
less than $\pi/2,$ and $s/t$ is sufficiently small, then
$d(\sigma_1(s),\sigma_2(t))<t,$ in particular $\sigma_1(0)=\sigma_2(0)$
does not minimize the distance from $\sigma_2(t)$ to the segment
$\sigma_1([0,\varepsilon])$.  We need this for the properties of
projection map below.  

The sum of angles of a geodesic triangle is less than or equal to $\pi$
( \cite{Ba}, I.5.2).  It follows immediately from this property that
if $I,J$ are geodesic segments issuing the same point and the subtended
angle equals $\pi$ then the concatenation of these segments is also a
geodesic segment.  If $\sigma$ is a geodesic segment (possible infinite)
then for any $x\in X$ there is (\cite{Ba}, I.5.6) a unique point
$p_\sigma(x)\in \sigma$ such that $d(x,p_\sigma(x))=d(x,\sigma)$.  The
map $p_\sigma$ is called the {\em projection} onto $\sigma.$ It follows
from the cosine formula and the remarks above that for each $x\in X$, the angles of
$[x,p_\sigma(x)]$ with $\sigma$ both are greater than or equal to $\pi/2$.

\subsection{Fake zero angles}
 We say that a complete $\CAT(0)$ space $X$ 
{\em has  fake zero angles} if there are two geodesics
issuing the same point, are disjoint except at that point, and 
the angle subtended at that point is zero. 

\begin{proposition}
\label{propsec}
A piecewise Euclidean $\CAT(0)$ complex $X$ with finitely many isometry
types of cells has no fake zero angles.
\end{proposition}
 
\proof
The assumptions imply that the path metric on $X$ is geodesic and
complete.  The angles can be defined in terms of link distance
\cite{BB}.  Namely, let $X$ be a piecewise Euclidean complex, $x\in
X$.  The {\em link} $\Lk_xA$ of the Euclidean cell $A$ is the set of
unit tangent vectors $\xi$ at $x$ such that a nontrivial line segment
with initial direction $\xi$ is contained in $A$.  We define the link  
$link(x)$ of $x\in X$ by $link(x)= \cup_{A\ni x}\Lk_xA,$ where the union is taken over
all closed cells containing $x$.  Angles in $\Lk_xA$ induce a natural
length metric $d_x$ on $link(x)$ which turns it into a piecewise
spherical complex.  The angle between $\xi,\eta\in link(x) $ is then
defined by $\angle_x(\xi,\eta)=\min(d_x(\xi,\eta),\pi).$ 

Any two segments
$\sigma_1,\sigma_2$ in $X$ with the same endpoint $x$ have the natural
projection image in the link of $x$ and $\angle_x(\sigma_1,\sigma_2)$
equals the angle between these two projections.  Now the assertion of
the lemma is clear since if the segments are disjoint, apart the origin,
then their images in the link are distinct and hence the link distance
is nonzero.
\endproof

\bigskip
\noindent
{\bf Example (V. Berestovskii): } 
Take $R^2$ with the positive $x$-axis removed and 
stick in the region $\{(x,y):x\geq 0,\ y\leq x^2\}$ along
the obvious isometry of the boundary. The result is a CAT(0)-
space with fake angles.

\section{Proof of Theorem \ref{theorem:main}}

We will need the following well-known lemma.

\begin{lemma}[Ping-Pong Lemma]
Let $\Gamma$ be a group of permutations  on a set $X$, let
$g_1,g_2$ be the elements of $\Gamma$ of order at least 3.
 If  $X_1,X_2$ are disjoint subsets of $X$ and for all $n\not= 0,\ i\not= j$,
$g_i^nX_j\subset X_i$
then $g_1,g_2$ freely generate the free group $F_2$.
\end{lemma}

The proof of Theorem \ref{theorem:main} divides into two cases,
depending on which of the two hypotheses is assumed.

\bigskip
\noindent
{\bf Assuming (1): }  Let $s$ be the midpoint of $S$.  Let $D_1,D_2$ be the
fundamental domains for the action of $g_1,g_2$ on $A_1,A_2$, 
chosen as open segments on
$A_1,A_2$ with the center $s.$ Let $p_1,p_2$ be the geodesic projection
maps of $X$ onto $A_1,A_2$ respectively.  Set $X_1=p_1^{-1}(A_1-D_2)$
and $X_2=p_2^{-1}(A_2-D_2).$ To apply the Ping-Pong Lemma we need
 to show that $X_1\cap X_2$ is empty.  Suppose, to the contrary,
it is not and let $x\in X_1\cap X_2.$  So $p_1(x)\in X_1, p_2(x)\in
X_2$.  Suppose that $p_1(x)\in A_1^+,p_2(x)\in A_2^-$, the other cases
being similar.  Consider the geodesic
triangle $xp_1(x)p_2(x)$.  The sum of its angles is at most $\pi$ and we
would like to get the contradiction with this.  
Suppose first that the angle at $x$ is nonzero.  
By angles assumption we get that $[p_1(x),p_2(x)]$ is the concatenation of
$[p_1(x),s_-],[s_-,s_+],[s_+,p_2(x)]$.  By the property of a projection
map the angles of the triangle at $p_1(x),p_2(x)$ are greater or equal
$\pi/2.$ Hence the sum of the angles is strictly greater than
$\pi$. This contradiction proves disjointness.
 
Suppose now that the angle at $x$ is zero, then by the fake zero angles
assumption, the geodesics $(x,p_1(x)],(x,p_2(x)]$ are not disjoint hence
$[x,p_1(x)]\cap [x,p_2(x)]=[x,y]$ for some $y\not=x.$ To get the
contradiction consider the triangle $yp_1(x)p_2(x)$.  If $y$ is strictly
closer to $x$ than both of $p_1(x),p_2(x)$, then the $y$-angle of the
triangle $yp_1(x)p_2(x)$ is nonzero and the sum of angles in this
triangle is strictly greater than $\pi$ - contradiction.  If not then
say $y=p_1(x).$ Since $y$ lies on $[x,p_2(x)]$, its projection onto
$A_2$ is the same as that of $x.$ Hence $p_2(x)=p_2(y).$ By definition,
$p_2(y)$ is the point on $A_2$ closest to $y.$ But $s_+\in A_2$ lies on
the geodesic $[y,p_2(x)]$ and $s_+\not=p_2(x)$; again, this is a
contradiction.

Finally it remains to check that $g_i^nX_j\subset X_i, i\not=
j, n\not= 0.$ Note first that $g_i$ commutes with $p_i$.  Indeed for any
$x\in X,$ $p_i(x)$ is the unique point in $A_i$ such that
$d(x,p_i(x))=d(x,A_i).$ But
$$d(g_ix,A_i)=d(g_ix,g_iA_i)=d(x,A_i)=d(x,p_i(x)) =d(g_ix,g_ip_i(x))$$
That is, the point $g_ip_i(x)$ realizes the distance $d(g_ix,A_i)$ and
thus it is the projection of $g_ix$ hence $g_ip_i(x)=p_ig_i(x)$.

Clearly $g_i^n D_i \subset A_i-D_i, n\not= 0$, whence
$g_i^n(X-X_i)\subset X_i$.  Finally, $X_j\subset X-X_i$, hence
$g_i^n(X_j)\subset X_i$.

\bigskip
\noindent 
{\bf Assuming (2):}
 Let $p_1,p_2$ be the geodesic projection maps of 
$X$ onto  $A_1,A_2$ respectively.
 Let $D_1,D_2$ be the fundamental domains for $g_1,g_2$ chosen
as the open segments on $A_1,A_2$ with the centers $a_1,a_2$
respectively.
 Let $p_1,p_2$ be the geodesic projection maps of 
$X$ onto  $A_1,A_2$ respectively.
 Set $X_1=p_1^{-1}(A_1-D_2)$ and $X_2=p_2^{-1}(A_2-D_2)$
and repeat the argument of the first case.
\endproof

\section{A family of examples}
\label{section:examples}

\subsection{Motivation: the Burau representation}
The braid group on $n$ strands, denoted $B_n$, is the group with
generators $s_1, s_2, \dots, s_{n-1}$ and relations
$$s_js_k=s_ks_j, s_is_{i+1}s_i=s_{i+1}s_is_{i+1}$$ for all possible $i,
j, k$ with $|k-j|\ge 2$ and $i+1,j,k\le n$.  The {\em Burau
representation} $\s_i\mapsto \s_i$ 
is a natural representation of $B_n$ on the $n$-dimensional linear space 
$V={K}^n$ over the field ${K}=\Q( t)$; in the standard basis $\{ e_i | i=1,
\dots ,n \}$ the Burau representation is determined by 
$$
\s_i(e_j)=
\left\{
\begin{array}{ll}
(1- t)e_i+e_{i+1}& \mbox{if\ }j=i\\
te_i& \mbox{if\ } j=i+1\\
e_j&\mbox{if\ }j\neq i,i+1
\end{array}\right.
$$

The subspace ${K}e$, where e=$\sum_{i=1}^n e_i$, is clearly an invariant
subspace for the group action, and the group acts trivially there.  Thus
the quotient space $V/{K}e$ has a natural action of the braid group.  
In the induced basis $e_1, e_2,\dots ,e_{n-1}$ the elements $\s_i$ 
act via matrices which we denote by $b_1, b_2,\dots ,b_{n-1}$.  
The only new aspect in this {\em reduced Burau representation} is that
$$b_{n-1}(e_{n-1}) = (1- t)e_{n-1} + e_n = -(\sum_{i=1}^{n-2} e_i +
te_{n-1}).$$ Our interest is in the case of four strands.  In this case
there is the following well known result (however, 
the matrices are incorrectly specified in \cite{Bir}).

\begin{proposition}[\cite{Bir}, Theorem 3.19]
The reduced representation of $B_4$ is faithful iff it is faithful
on the free group generated by $a=s_3s_1^{-1}$ and
$b=s_2(s_3s_1^{-1})s_2^{-1}$.
\end{proposition}

Let $f$ (resp. $k$) be the image of $a$ (resp. $b$) in $\GL_3(K)$ under
the reduced Burau representation.  It is not difficult to see that both
$f$ and $k$ are diagonalizable.  In fact, by conjugating the Burau
representation, and changing $t$ to $-t$, we may take $f$ to be the
diagonal matrix with entries $1,- t^{-1}, - t$, and $k$ becomes the
matrix $k=sfs^{-1}$ where $$s=(1- t)^{-2}
\left(
\begin{array}{ccc}
-(1+ t)&1+ t^2&- t(1+ t^2)\\
1&- t& t\\
1&-1&t^2
\end{array}
\right)
$$

We consider the action of $\GL_3({K})$ on the Bruhat-Tits building
${\cal T_\val}$ for the field ${K}$ with the discrete valuation
$\val=\val_{\infty}$ at infinity.  We have two elements of $\GL_3({K})$
which are acting so that they each stabilize an entire apartment of
${\cal T}$; these apartments $A_u$ and $A_v$ on general principles will
meet in a convex subset of the building.

However, by analogy with the case of actions on
trees, we might expect that if the intersection of $P_f$ and
$P_k$ is sufficiently small with respect to the translation distances of
$f$ and $k$ then the group generated by $f$ and $k$ is free.  Since $f$
is semisimple it is
easy to see that it acts by translation  on its apartment not along any
wall and
similarly for $k$.  We can determine  the intersection 
$P_f\cap P_k$.

\begin{lemma}[$P_f\cap P_k$ is a point]
\label{lemma:point}
The intersection of $P_f$ and $P_k$ consists of
precisely one lattice class.
\end{lemma}

\proof Let $\val$ denote the discrete valuation, with 
valuation ring $\cal O$ and uniformizer $\pi$.  
The lattices which represent the lattice classes in $P_f$ are
$$L_{a_1,a_2,a_3}={\cal O}\pi^{a_1}e_1 + {\cal O}\pi^{a_2}e_2 + {\cal
O}\pi^{a_3}e_3.$$ 
Since $k=sfs^{-1}$, the lattices in $P_k$ are precisely
the lattice classes $[sL]$ for $[L]$ in $P_f$, the standard apartment. 
Thus
a lattice class of $L_{-a_1,-a_2,-a_3}$ belongs to the intersection 
if there are integers $a_1, a_2, a_3$ and $b_1, b_2, b_3$ so
that $[sL_{b_1,b_2,b_3}]=[L_{-a_1,-a_2,-a_3}]$.  In other words
there is  some matrix $m$ in $\GL_3({\cal O})$ so that 
$$s
\left(
\begin{array}{ccc}
\pi^{b_1}&0&0\\
0&\pi^{b_2}&0\\
0&0&\pi^{b_3}
\end{array}\right)
=
\left(
\begin{array}{ccc}
\pi^{-a_1}&0&0\\
0&\pi^{-a_2}&0\\
0&0&\pi^{-a_3}
\end{array}
\right)
m$$  
This gives rise to the conditions: 
$\pi^{a_i}s_{i,j}\pi^{b_j}\in{\cal O}$ for all entries
$s_{i,j}, 1\leq i,j\leq 3$ of $s$.  This implies 
$$a_i + b_j + \val(s_{i,j}) \ge 0$$

Since $m$ is invertible, $\val(\det(m))=0$.  This implies
$$a_1 + a_2 + a_3 + b_1 + b_2 + b_3 + \val det(s)=0.$$ 
We can convert many of the above inequalities into equalities by the 
following argument.  Direct calculation gives $\val(\det(s))=2$.  Hence 
$$2=(a_2+b_1)+(a_1+b_2)+(a_3+b_3)\ge 2+0+0$$ and therefore $a_2+b_1=2$,
$a_1+b_2=0$, and $a_3+b_3=0$.  It follows similarly that $a_1+b_1=1$,
$a_2+b_2=1$, $a_3+b_1=2$, $a_1+b_3=-1$, $a_2+b_3=0$, $a_3+b_2=1$.  From
this we
can immediately solve the nine equations in six variables to get the
solutions
$a_1=-1-c, a_2=-c, a_3=-c, b_1=2+c, b_2=1+c, b_3=c$ and so this
determines
exactly one lattice class solution to the intersection of the two
apartments.
\endproof

\subsection{A $2$-parameter family of representations}

Based on the example afforded by the Burau representation, we consider
the transformations $f$ and its conjugate $k=sfs^{-1}$, where  
$$s=(1- t)^{-2}
\left(
\begin{array}{ccc}
-(1+ t)&1+ t^2&- t(1+ t^2)\\
1&- t& t+\beta t^2\\
1&-1+\alpha t& t^2\\
\end{array}
\right)
$$ for $\alpha$ and $\beta$
 any rational numbers. We can think of this as giving a $2$-parameter
family of representations 
$$\rho_{\alpha,\beta}:F_2\rightarrow
\GL_3(\Q(x))$$
\noindent
for $\alpha,\beta\in\Q$.  Consider any fixed parameters, giving a 
pair $f,k$.  As in the proof of Lemma 
\ref{lemma:point}, the invariant apartments of $f$ and $k$ 
meet at exactly one point.

The lattices which represent the lattice classes in the standard
apartment $A$
are
$L_{a_1, a_2, a_3}={\cal O}\pi^{a_1}e_1 + {\cal O}\pi^{a_2}e_2 + {\cal
O}\pi^{a_3}e_3,$ $\pi$ represents the
uniformizer
$ t^{-1}$. This apartment is stabilized by $f$. In this standard basis,
$e_1, e_2, e_3$,  $f$ is represented by a diagonal matrix with diagonal
entries $1,  t^{-1},  t$. The apartment
stabilized by
$k=sfs^{-1}$ is $sA$. We consider the lattice class $x$ of the lattice
$L_{-1,0,0}$ in the standard apartment as the special vertex which is
the
common cone point of our sectors
$P_f^{+},P_f^{-},P_k^{+},P_k^{-}$.  This vertex is also
$sL_{-2,-1,0}$ by the calculation above. The walls of the sectors
$P_f^{+}$
and $P_f^{-}$ are then easily seen
to be the subcomplex represented by the lattice classes of
$M_n=L_{-1, n, 0}\ n\in Z$, and $N_n=L_{-1, 0, n}\ n\in Z$. The
translation
$f$ takes $[L_{-1, 0, 0}]$ to
$[L_{-1, 1, -1}]$ with axis bisecting the walls of the sector. We can
apply
$s$ to this configuration to
obtain a similar description for the axis and sectors for the element
$k$.

We shall calculate the link of the vertex
$x$ and represent it in terms of the spherical building of $GL_3(\bf
Q)$,
where the rational field $\bf Q$
is the residue field of the  field
${\bf Q}( t)$ with respect to the valuation at infinity. A $q$-simplex
in the building is represented by
a chain of lattices $L_0\subset L_1\subset\cdots\subset L_q$ with $\pi
L_q\subset L_0$ and
$L_{r+1}/L_r\cong {\bf Q}, \ r\ge 0$. The link of the vertex $[L]$ is
the
simplicial complex  whose
vertices are lattice classes $[L^{'}]$ so that
$\pi L\subset L^{'}\subset L$. By taking chains of such classes we
obtain
a simplex in
the spherical building of
$GL_3(\bf Q)$, viewed as flags in the 3-dimensional ${\bf Q}$-vector
space, $L/\pi L$.

With the labelling of the walls as above, the first chamber in the
sector
$P_f^{+}$  is represented by
the lattice classes of $L=L_{-1, 0, 0}$ and the lattice classes of
$M_{1}=L_{-1, 1, 0}$ and
$N_{-1}=L_{-1, 0, -1}$ while the first chamber in $P_f^{-}$ is
represented by
the lattice classes of
$L_{-1, 0, 0}$ and the lattice classes of
$M_{-1}=L_{-1, -1, 0}$ and $N_{1}=L_{-1, 0, 1}$. We can apply $s$ to
these
classes to represent the other
leading chambers in $P_k^{+}$ and $P_k^{-}$.

In order to obtain the precise configurations in $P_f^{+}$, considering
the
leading chamber  which contains the lattice classes of $M_1$ and
$N_{-1}$. We have
the chain of lattices $\pi L\subset \pi N_{-1}\subset M_1\subset L$, so
in $link(x)$ we have the flag of subspaces

$$X_1=\pi N_{-1}/\pi L\subset X_2= M_1/\pi L,$$
which is the line $X_1$ with basis $e_3$ as a
subspace of $X_2$ with basis $\{e_1,e_3\}$. The same can be done for the
leading chamber of $P_f^{-}$ to
obtain the flag
$$\pi L\subset\pi M_{-1}\subset N_1\subset L$$ giving the line $Y_1$
with
basis $e_2$ as a subspace of
$Y_2$ with basis $\{e_1,e_2\}$. It is immediate that the edge
$X_1\subset
X_2$ is opposite to $Y_1\subset Y_2$.

Similarly, we have the flags $s(X_1)\subset s(X_2)$ and $s(Y_1)\subset
s(Y_2)$ in $link(x)$ coming from the sectors $P_k^{+}$ and $P_k^{-}$
based at
the point $[s(L_{-2,-1,0})]=t$,
since $[L_{-1,0,0}]=[s(L_{-2,-1,0})]$. Consider the change of basis
matrix described in the proof of Lemma \ref{lemma:point},  
$$m=(1- t)^{-2}
\left(
\begin{array}{ccc}
-t(1+t)&1+t^2&-(1+t^2)\\
t^2&-t^2&t\\
t^2&-t&t^2\\
\end{array}
\right)$$ 
and reduce mod $\pi=\frac{1}{t}$ to obtain
the rational matrix
$$\left(
\begin{array}{ccc}
-1&1&-1\\
1&-1&\beta \\
1&\alpha  &1
\end{array}\right)
$$
to find that $s(X_1)$ has basis $(-1,\beta,1)$ and $s(X_2)$ has basis
$\{(-1,1,1),(-1,\beta,1)\}$
and $s(Y_1)$ has basis $(1,-1,\alpha)$ while $s(Y_2)$ has basis
$\{(1,-1,\alpha),
(-1,\beta,1) \}$. Thus these edges are {\it opposite}
to the edge $X_1\subset X_2$ for $\alpha$ and $\beta$, rationals which
are neither 0 nor 1 ; but {\it not opposite} if $\alpha$ or $\beta$ is 0
or 1. 

To see this we can describe the local hexagons in the link. Oppositeness
in this case means that the 
two 2-dimensional subspaces intersect along a line which is not the
special line in each.  We  show that $X_1\subset
X_2$ is opposite to $s(X_1)\subset s(X_2)$ and $s(Y_1)\subset
s(Y_2)$. For example, $X_2$ is spanned by $e_3$ and
$e_1$, and $sX_2$ is spanned by $(-1,\beta,1)$ and  $(-1,1,1)$. The
subspace $X_2\cap s(X_2)$ does not contain either the line
generated by $e_3$ or the line generated by $(-1,\beta,1)$ iff
$\beta\ne0$. Thus we have oppositeness in this case. Similarly, we can
treat the cases of $X_2$, $s(Y_2)$, and $Y_2, s(Y_2)$.

\bigskip
\noindent
Roger C. Alperin:\\
Dept. of Mathematics and Computer Science, San Jose State University\\
1 Washington Square \\
San Jose, CA 95192\\
E-mail: alperin@mathcs.sjsu.edu
\medskip

\noindent
Benson Farb:\\
Dept. of Mathematics, University of Chicago\\
5734 University Ave.\\
Chicago, IL 60637\\
E-mail: farb@math.uchicago.edu
\medskip

\noindent
Guennadi A. Noskov:\\
Institute of  Mathematics, Russian Academy of Sciences\\
Pevtsova 13\\
Omsk 644099, RUSSIA\\
E-mail: noskov@private.omsk.su

\end{document}